\documentstyle[12pt]{amsart}

\input xy
\xyoption{all}

\newtheorem{thm}{Theorem}[section]

\newtheorem{propn}[thm]{Proposition}
\newtheorem{rem}[thm]{Remark}
\newtheorem{lemma}[thm]{Lemma}
\newtheorem{cor}[thm]{Corollary}

\newcommand{\Hom}{\operatorname{Hom}}

\renewcommand\Im{\operatorname{Im}}
\newcommand\Id{\operatorname{Id}}
\newcommand\Pic{\operatorname{Pic}}

\newcommand\cw{{\mathcal W}}
\newcommand\cp{{\mathcal P}}
\newcommand\cm{{\mathcal M}}
\newcommand\cq{{\mathcal Q}}
\newcommand\co{{\mathcal O}}
\newcommand\cu{{\mathcal U}}
\newcommand{\N}{{\mathcal {N}}}

\begin{document}

\title [On Picard Bundles over Prym varieties] {On Picard Bundles over
Prym varieties}

\author{L.~Brambila-Paz}

\address{Leticia Brambila-Paz. CIMAT, Apdo. Postal 402, C.P. 36240
Guanajuato. M\'exico.}

\email{lebp@@fractal.cimat.mx}

\thanks{The authors are members of the VBAC Research
group of Europroj. The first author
acknowledges support from CONACYT (Grant no. 28-492-E).
The second author acknokwledges
support from the Spanish DGES through the research
project  PB96-1305 and by the ``Junta de Castilla y Le\'on'' through the
research
project SA27/98. The third author acknowledges support from CIMAT 
Guanajuato, GNSAGA Italian Research Council}

\author{E.~G\'omez-Gonz\'alez}

\address{Esteban G\'omez-Gonz\'alez. Departamento de Matem\'aticas.
Universidad de Salamanca. Plaza
de la Merced 1-4. 37008 Salamanca. Spain}

\email{esteban@@gugu.usal.es}

\author{F.~Pioli}

\address{Fabio Pioli. Dipartimento di Matematica, Universit\`a degli Studi
di Genova, Via
Dodecaneso 35, 16146 Genova, Italy}

\email{pioli@@dima.unige.it}



\maketitle

\begin{abstract}
Let $P$ be the Prym variety associated with a covering $\pi : Y \to X$ between
non-singular irreducible projective curves. If $\widetilde P$ is a principally
polarized Prym-Tyurin variety associated with $P$, we prove that the induced 
Abel-Prym morphism $\widetilde \rho \colon Y\to \widetilde P$ is birational 
onto its image for genus $g_X>2$ and $\operatorname{deg} \pi\neq 2$. We
use such result to prove that the Picard bundle over the Prym variety is
simple and moreover is stable when $\widetilde \rho $ is not birational. As a
consequence we obtain that the Picard bundle on the moduli
space $\cm_X (n,\xi)$ of stable vector bundles with fixed determinant and rank 
$n$ is simple for $g_X\geq 2$ and in the case $g_X = 2$ and $n=2$ then we 
prove that the Picard bundle on $\cm_X (2,\xi)$ is stable. 
\end{abstract}


\section*{Introduction}
Let $X$ be a non-singular irreducible projective curve of genus
$g_X\geq 2$. We denote by $\cm (n,d)$
the moduli space of stable vector bundles over $X$ of rank $n$ and degree
$d$. If $n$ and $d$ are coprime there exists a
universal family $\cu $
parametrized by $\cm (n,d)$, called the Poincar\'e bundle. The higher direct
images of the Poincar\'e bundle on $\cm (n,d)$
are called Picard sheaves. If $d> 2n(g_X-1)$, the direct
image $\cw$ of $\cu$ is locally free and it is called the Picard bundle.

For $n=1$, $\cm (1,d)$ is the Jacobian $J^d(X)$. Ein and Lazarsfeld in
\cite{EL} proved the stability of the Picard bundle $\cw _J$
 when $d>2g_X-1$ and
Kempf (see \cite {Ke}) for $d=2g_X-1$.
For $n\geq 2$, Li  proved  the stability of $\cw $ over $\cm (n,d)$ when
 $d > 2g_Xn$ (see \cite{Li}).
Denote by $\cw _{\xi}$ the
 restriction of $\cw$ to the subvariety $\cm _{\xi}\subset \cm (n,d)$
determined by stable bundles with fixed determinant $\xi \in J^d(X)$.
Balaji and Vishwanath  proved in \cite{bv} that for rank 2, $\cw _{\xi}$ is
simple. Actually, they compute the deformations of $\cw _{\xi}$.

In this paper we prove (see Theorem 3.1) that for any rank
$\cw _{\xi}$ is simple
when $g_X\geq 2$ and $d>n(n+1)(g_X-1)+6$ and in the case $g_X=2$ and $n=2$ then 
we prove that
$\cw_\xi$ is  stable for $d>10$. 

The result follows from studying the restriction of the Picard bundle
over the Jacobian of a curve (which is a covering space of $X$) to some natural 
subvarieties,namely the Prym variety and the Jacobian of the base space. That is, 
let
$\pi :Y\rightarrow X$ be a covering between  non singular irreducible
projective curves. Let $(J(Y), \Theta_Y)$ be the principally polarized Jacobian 
of degree
zero of $Y$ and $\cp _Y$ the Poincar\'e bundle over $Y\times J(Y)$. Fixing
a line
bundle $L_0$ on $Y$ of degree $d>2g_Y-2$, we identify $J(Y)$ with $J^d (Y)$
and consider the Picard bundle $\cw_J =p_{2*}(p_1^*L_0\otimes
\cp _Y)$ on $J(Y)$. We denote the restriction of $\cw_J$ to the subvarieties
$P$ and $\pi^* (J(X)) $ of $J(Y)$ by $\cw_P$ and $\cw_\pi$ respectively.

 We prove  that if $\pi_* (L_0)$ is stable,  then $(\pi^*)^*(\cw_J)$ is
$\Theta_X$-stable
on $J(X)$ where  $\pi^*\colon J(X)\to J(Y)$ is the pull-back morphism (see
Theorem
\ref{th:beau} and Remark \ref{rem:beau}). We deduce that
$\cw_\pi$ is stable with respect to the polarization
$\Theta_Y |_{\pi^* (J(X))}$.

The restriction $\Theta _P$ of $\Theta _Y$ to $P$ need not
be a multiple of a principal polarization. However, it
is possible to construct (see \cite{We}) a
principally polarized abelian
variety $(\widetilde P, \varXi)$ and an isogeny $f\colon \widetilde P \to
P$ with $f^{-1}
(\Theta_P)\equiv n\varXi$ such that there exists a map $\widetilde \rho
\colon Y\to \widetilde P$ with $\widetilde \rho (Y)$ and $\varXi^{m-1}$
numerically equivalent where $m=\dim P = \dim \widetilde P$.

We prove that if $n\neq 2$ and $g_X>2$, the morphism
$\widetilde \rho \colon Y\to
\widetilde P$ is birational onto its image
(Theorem \ref{th:bir})
 and  we use some properties of the Fourier-Mukai transform to prove
that in this case $\cw _P$ is simple.

Moreover, we give explicitly the cases where $\widetilde \rho$ may not be
birational. In this case, we prove that the Prym variety is the image in
$J(Y)$ of the Jacobian of the normalization of the curve $\widetilde \rho
(Y)$. As a consequence of stability of $\cw_\pi$, we have that the Picard bundle
$\cw_P$ is stable with  respect to $\Theta_Y |_P$ when $\widetilde \rho $ is not
birational. In particular, for $g_X\geq 2$ and $d> 2g_Y+4$,
$\cw_P$ is simple (Theorem \ref{th:sim}).

{\small {\it Acknowledgments} : We thank Laura Hidalgo for previous 
discussions and G.P. Pirola for useful suggestions concerning Theorem
\ref{th:bir}. The first and the third named authors acknowledge the generous 
hospitality of the Departamento de Matem\'aticas of the University of 
Salamanca. The third author wishes to thank the Centro de Investigaci\'on en 
Matem\'aticas
(CIMAT), Guanajuato for the warm hospitality and support.}

\section{Prym varieties}\label{sec:uno}

We shall denote by $J(X)$ the Jacobian of degree zero
of a non-singular projective irreducible curve $X$  and by $\Theta_X$ 
the natural polarization on $J(X)$ given by the Riemann theta divisor.

Let $\pi\colon Y\to X$ be a covering of degree $n$ between non-singular
irreducible projective curves of
genus $g_Y $ and $g_X$ respectively. We have the norm map $\operatorname{Nm}\pi \colon  J(Y) \to J(X)$ and the pull back map $\pi^* \colon J(X)
\to J(Y)$. 
It is known that the map $h \colon J(X)\to \pi^* (J(X))$ induced by $\pi^*$ is an
isogeny such that
$h^{-1} \left( \Theta_Y |_{\pi^* (J(X))} \right) \equiv n
\Theta_X$. We will denote by $\Theta _{\pi}$
the restriction of the theta divisor of $J(Y)$ to
$\pi^*(J(X)).$

The Prym variety $P$ associated with the covering is the
abelian subvariety of $J(Y)$
defined by
$P= \Im (n\Id_{J(Y)} - \pi^*\circ \operatorname{Nm}\pi)$.
We shall denote by $\mu $ the morphism $(n\Id_{J(Y)} - \pi^*\circ
\operatorname{Nm}\pi)\colon J(Y)\longrightarrow P$.

On $P$ there exists a natural polarization $\Theta _P$ given by the
restriction of the theta
divisor $\Theta_Y $ to $P$. In general $\Theta _P$ is not a multiple of
principal polarization (see
\cite{LB}, Theorem 3.3, pp. 376). The following theorem gathers several
standard results about Prym-Tyurin varieties for which we refer to
\cite{LB}, \cite{We} and \cite{Go}.

\begin{thm}\label{th:isogenia} Let $\pi \colon Y\to X$ be a covering of degree
$n$. There exists a principally polarized abelian variety $(\widetilde P,
\varXi) $,
an isogeny $f:\widetilde P\to P$ and a morphism $\widetilde \mu \colon J(Y)
\to \widetilde P$ such that: \newline
\indent {\rm 1.} $f^{-1} (\Theta_Y |_P) \equiv n \varXi $.\newline
\indent {\rm 2.} The map $f\circ \widetilde \mu \colon J(Y) \to P$ coincides
with the morphism
$\mu$ defined as above.\newline
\indent {\rm 3.} The restriction $\widetilde \mu |_{P}$ is surjective and the
composition map $f\circ
\widetilde \mu |_{P}$ coincides with the multiplication by $n$ on $P$.

Moreover, the following relation holds
$$ \operatorname{Nm}^{-1}\pi (\Theta_X) + {\widetilde \mu }^{ -1} (\varXi )  
\equiv  n \Theta_Y,$$ where $\equiv $ denotes algebraic equivalence between 
divisors.
\end{thm}

 From now on we shall consider a fixed isogeny $f\colon  \widetilde P \to P$
enjoying  the properties listed in Theorem 1.1. Let $y_0$ be a fixed point
on $Y$. Let
${\widetilde \rho}_{y_0} \colon Y \to \widetilde P $ and
 $\rho_{y_0} \colon Y \to P$ be the morphisms
$${\widetilde \rho}_{y_0} = \widetilde \mu \circ \alpha_{y_0} \qquad \qquad
\text{and} \qquad
\qquad \rho_{y_0} = f\circ {\widetilde \rho}_{y_0} =
\mu \circ \alpha_{y_0}$$
where $\alpha_{y_0}\colon Y \to J(Y)$ is the Abel-Jacobi map defined by
$y\mapsto \co (y-y_0)$.
We shall refer to $\rho_{y_0} $ and
$\widetilde {\rho}_{y_0} $ as {\it the Abel-Prym maps} of
$P$ and $\widetilde P $ respectively. When there is no confusion on the
point $y_0$ we shall
omit it in the notations.

By the theory of Prym-Tyurin varieties (see \cite{We} and \cite{LB}) the
morphism
${\widetilde \rho}$ has the following property:
$${\widetilde \rho}_{*} [Y] =
\frac n{(m-1)!} \bigwedge^{m-1} [\varXi],$$
where $m=\dim P =\dim \widetilde P $. This last property is usually
expressed by
saying that the curve $Y$ (or rather the morphism
${\widetilde \rho}$) is of class $n$ in $(\widetilde P, \varXi)$.

In general the Abel-Prym map $\rho \colon Y \to P$ need not be birational.
However, we have the following theorem.

\begin{thm}\label{th:bir} Let $\pi \colon Y\to X$ be a covering of degree
$n$. If $g_X>2$ and
$n\neq 2$ then Abel-Prym
morphism
${\widetilde \rho} \colon Y\to \widetilde P$ is birational onto its
image.
\end{thm}

In order to prove Theorem \ref{th:bir} we need a lemma.
Let $(A, \Theta_A)$ be a principally polarized abelian variety of dimension
$d$ and $C$ a non-singular
projective irreducible curve. Let $\nu \colon C\to A$ be a morphism and let
$u\colon J(C) \to A$
be the unique morphism (up to a translation) obtained by the universal
property of the Jacobian such
that $\nu = u \circ \alpha$. Using the polarizations, we define the
transposed morphism
$u^t\colon A\to J(C)$ by
$$ u^t = \phi_{\Theta_C}^{-1} \circ \hat u \circ \phi_{\Theta_A},$$
where $\hat {u}\colon \hat {A} \to \widehat{J(C)}$ is the dual morphism of
$u$ and
$\phi_{\Theta_A}\colon A\to \hat A$, $\phi_{\Theta_C} \colon J(C) \to
\widehat{J(C)}$ are the
isomorphisms induced by the principal polarizations $\Theta_A $ and
$\Theta_C$ respectively.

Recall that a curve $C$ is of class $n$ in $(A, \Theta_A)$, that is $\nu_* [C]
= \frac{n}{(d-1)!}\bigwedge^{d-1} [\Theta_A]$, if and only if $ u\circ u^t
= n \Id_A$ (see
\cite{We}).

\begin{lemma}\label{lem:div} Let $(A,\Theta_A)$ be a principally polarized
abelian variety and $\nu \colon C\to A$ a curve of class $n$ in $A$. Let
$Z$ be the normalization of the image $\nu (C) $ and let $p\colon C\to Z$
and $\widetilde \nu \colon Z\to A$ be the induced maps. If $r$ is the
degree of $p$, then $r$ divides $n$ and the curve $Z$ is of class $\frac
{n}{r}$
in $A$.
\end{lemma}

\begin{pf}
Let $\operatorname{Nm} p\colon J(C)\to  J(Z)$ be the norm map of $p$ and $u\colon 
J(C)\to
A$ the morphism
corresponding to $\nu $ obtained by the universal property of the Jacobian.
Similarly, $\widetilde u \colon J(Z)\to A$ corresponds to $\widetilde \nu$.
  Since $\nu \colon C\to A$ is of class $n$ and
$N_p\circ p^*= N_p\circ N_p^t = r
\Id_{J(Z)}$ ,
we have:
$$ n\Id_A = u \circ u^t = \widetilde u \circ  N_p \circ N_p^t\circ
{\widetilde u}^{t}= r\cdot \left(
\widetilde u \circ {\widetilde u}^t \right).$$

Thus, for every $x\in A_r$ in the $r$-torsion of $A$,  $r\cdot (\widetilde u
\circ {\widetilde u}^t) (x)
= 0$ and therefore $nx=0$. Hence, $A_r\subset A_n$ and $r$ divides $n$.
Since $A$ is connected and
$A_r$ is discrete, the morphism $A\to A_r$ given by $x\mapsto \left(
\widetilde u \circ
{\widetilde u}^t \right) (x) - \frac {n}{r} x $ is identically zero. Hence
$\widetilde u \circ {\widetilde u}^t = \frac{n}{r}\Id_A$.
\end{pf}

{\it  Proof of Theorem \ref{th:bir}.}  Suppose that
$\widetilde {\rho}$ is not birational onto
its image and $n\neq 2$.
Let  $Z$ be the normalization of the curve $\widetilde {\rho}
(Y)$ and let
$\widetilde \pi \colon Y\to Z$ be the morphism induced by ${\widetilde
\rho}$. Let
$\tilde n$ be the degree of $\tilde \pi$. That is,
\begin{equation*}
\xymatrix{ & Y\ar[dl]_{n:1}^{\pi}
\ar[dr]_{\tilde
\pi}^{\tilde n:1} &
\\ X & & Z}
\end{equation*}
Since the morphism $\tilde \mu$
is surjective, $Z$
generates $\tilde P$. Therefore the map $\tilde u \colon J(Z)\to \widetilde
P$ defined by the universal
property of the Jacobian is surjective. Hence $g_Y-g_X = \dim \widetilde
P \leq g_Z$, where
$g_Z$ is the genus of $Z$.

 By Riemann-Hurwitz formula for the coverings
$\pi$ and $\tilde \pi$, we have the following inequality
\begin{equation}\label{eq:genus}
 g_Y\leq \frac {g_Y-1-(\operatorname{deg} R_\pi)/2} { n }  +
\frac { g_Y-1- (\operatorname{deg} R_{\tilde \pi})/2 }
{\tilde n}  + 2
\end{equation}
where $R_\pi$ and
$R_{\tilde \pi}$ are the ramification divisors.

{}From Lemma \ref{lem:div}, it follows that $\tilde n$
divides $n$; let $q=n/\tilde n$. Hence, inequality (\ref {eq:genus})
implies that
$$
(n-q-1) g_Y\leq 2n -q-1 -\left( (\operatorname{deg} R_\pi)/2 + q
(\operatorname{deg} R_{\tilde \pi}) /2 \right).
$$
Since $\left( (\operatorname{deg} R_\pi)/2 + q (\operatorname{deg}
R_{\tilde \pi}) /2 \right) \geq
0$ and $n\neq q+1$ because $n\neq 2$, we have
 $$ g_Y\leq 2 + \left[ \frac{q+1}{n-q-1}\right]$$
where $[-]$ denotes the integral part of a rational number. Since $\tilde
n>1$, it is easy to
check that $\left[ \frac{q+1}{n-q-1}\right]\geq 1$ if and only if $\tilde
n=2,3,4.$ Studying
these cases we obtain that if ${\widetilde \rho}$ is not birational
and $n\neq 2$, then
\begin{alignat*}{2}
n&=3,\quad \tilde n=3 \text{ and } g_Y \leq 4; \qquad & n=4,
&\begin{cases} \text{ if }
\tilde n=4 &\text{ then } g_Y \leq 3\\
\text{ if } \tilde n = 2 &\text{ then } g_Y \leq 5
                \end{cases};\\
n&=6 , \begin{cases} \text{ if } \tilde n=6 &\text{ then } g_Y\leq 2\\
   \text{ if }  \tilde n=3 &\text{ then } g_Y \leq 3\\
    \text{ if } \tilde n=2 &\text{ then } g_Y\leq 4
     \end{cases};
\qquad\qquad &  n\neq 3,4,6,  &\begin{cases}
\text{ if } \tilde n=2 &\text{then } g_Y \leq 3\\
\text{ if } \tilde n\neq 2 &\text{then } g_Y \leq 2
                \end{cases} .
\end{alignat*}

In all cases, we have that if ${\widetilde \rho}$ in not birational
onto its image and
$n\neq 2$, then $g_Y\leq 5$.

Now if we also suppose that $g_X\geq 2$, by the Riemann-Hurwitz formula
for the covering $\pi$, we
deduce  that if $\operatorname{deg} \pi \neq 2$ and $g_X\geq 2$, then
${\widetilde \rho}$ is
not birational onto its image if the following cases occur:
\begin{equation}\label{eq:specials}
\begin{aligned}
&n=3, \quad \tilde n=3, \quad g_Y =4, \quad g_X=2, \quad g_Z=2, \quad
R_\pi=R_{\tilde \pi}=0
 \\   &n=4,  \quad \tilde n=2, \quad g_Y =5, \quad g_X=2, \quad g_Z=3,
\quad
R_\pi=R_{\tilde \pi}=0
\end{aligned}
\end{equation}
Therefore, if we assume that $g_X>2$ and $n\not= 2,$
the theorem is proved.
\hfill $\square$

\begin{rem}\label{rem:iper}\begin{em} From the proof of Theorem \ref{th:bir}
we deduce that if
$g_X\geq 2$, then the Abel-Prym
map
$\widetilde {\rho}$ may not be birational onto its
image possibly in the cases given in
 (\ref{eq:specials}) and when $n=2$. For $n=2$, it is well-known that if $Y$ is
not an hyperelliptic
curve then $\rho$ is birational onto its image and therefore
$\widetilde {\rho}$ is birational onto its image.
\end{em}\end{rem}

We finish this section with a lemma that we will use later.

\begin{lemma}\label{lem:lemcomp} For any $\tilde \xi \in \widetilde P $, we
have
$${\widetilde \rho}^* \left(\tau_{\tilde \xi}^* \co_{\widetilde P} (\varXi)
\right) \cong {\widetilde \rho}^* \co_{\widetilde P} (\varXi) \otimes
f(\widetilde \xi)^{-1},$$
where $\tau_{\tilde \xi}$ is the translation by $\tilde \xi$.
\end{lemma}

\begin{pf} By Theorem \ref{th:isogenia}, there exists $\xi \in P$ such that
$\tilde \mu (\xi)=
\tilde \xi$. Using the relation between the polarizations given in this
theorem, we obtain
 $$\begin{aligned}
{\widetilde \rho}^* \left(\tau_{\widetilde \xi}^* \co_{\widetilde P}
(\varXi) \right) &\cong \alpha^* \left({\widetilde \mu}^*
\left(\tau_{\widetilde \xi}^* \co_{\widetilde P} (\varXi)\right)\right)
\cong \alpha^* \left(\tau_{\xi}^* \left({\widetilde \mu}^* \co_{\widetilde P}
(\varXi)\right)\right)\\
&\cong \alpha^* \left(\tau_{\xi}^* \left(\co_{J(Y)} (n\Theta_Y)
\otimes N_\pi^* \co_{J(X)} (\Theta_X)^\vee \otimes \N\right)\right),
\end{aligned}
$$
where $\N \in \Pic^0 (J(Y))$. In particular, $\N$ is
invariant under translations. Moreover, we have $\tau_{\xi}^* \left(
N_\pi^* \co_{J(X)}
(\Theta_X)^\vee \right) \cong  N_\pi^*  \left(\co_{J(X)} (\Theta_X)^\vee
\right) $ because
$\xi \in P$ and the isomorphism
$$\alpha^* \left(\tau_{\xi}^* \co_{J(Y)} (n\Theta_Y) \right)\cong \alpha^*
\left(\co_{J(Y)} (n\Theta_Y) \right) \otimes \xi^{-n} .$$
Therefore, it follows that $\widetilde {\rho}^* \left(\tau_{\widetilde
\xi}^* \co_{\widetilde P}
(\varXi)\right) \cong \widetilde {\rho}^*\co_{\widetilde P} (\varXi)\otimes
\xi^{-n}$.
Since $\xi^{-n} = f (\widetilde{\mu} |_{P} (\xi^{-1})) = f(\tilde \xi)^{-1}$,
the lemma is proved.
\end{pf}

\section{Stability and restrictions of Picard bundles}\label{sec:dos}

Let $Y$ be a non singular irreducible projective curve of genus $g_Y\geq 2$. 
We fix a line bundle $L_0$ on $Y$ of degree
$d$ and a point $y_0\in Y$. Let $\cp_Y$ be the
Poincar\'e bundle over
$Y\times J(Y)$ normalized with respect to the point $y_0$,
i.e. $\cp_Y|_{\{y_0\}\times J(Y)}
\cong \co_{J(Y)}$.
The Picard sheaves (relative to $L_0$) on $J(Y)$ are defined as the higher
direct images  $R^ip_{2*}(p_1^* L_0 \otimes \cp_Y)$ where
$p_j$  are the canonical projections of $Y\times
J(Y)$ in the $j=1,2$ factor.
If $d>2g_Y-2$, then $\cw_J: = p_{2*}(p_1^* L_0 \otimes \cp_Y)$ is a vector bundle 
known as {\it
the Picard bundle}. We shall consider the restriction of 
$\cw_J$ to some subvarieties when the curve is a covering space.

 Let $\pi \colon  Y\to X$ be a covering of degree $n$ between non-singular
projective irreducible
curves of genus $g_Y$ and $g_X$ respectively and let
$\pi ^*\colon J(X)\rightarrow J(Y)$ be the pull-back morphism. As in \S 1
consider the Prym variety $(P,\Theta _P)$ and $ (\pi^*(J(X)),\Theta
_{\pi})$ in $J(Y)$.

{}From now on we fix a line bundle
$L_0$ on $Y$ of degree $d>2g_Y-2$.
Denote by $\cw _P$ and $\cw _{\pi}$ the restrictions of the Picard
bundle $\cw_J$ (relative to $L_0$) to $P$ and
$\pi^*(J(X))$ respectively.

\begin{rem}\begin{em}\label{rem:beau}  We can take a line bundle $L_0\in
J^d(Y)$
such that $\pi_*(L_0)$  is stable and
$\operatorname{deg} \pi_*(L_0) > 2ng_X$ without loss of generality.
Indeed: for a generic line bundle $L$ on $Y$,
A. Beauville in \cite{Be} has proved  that $\pi_* L$ is
stable on $X$ if
 $|\chi (L)|\leq g_X+\frac{g_X^2}n$
or $\operatorname{deg}\pi < \operatorname{max}\{ g_X(1+\sqrt 3)-1,2g_X+2\}$.
Therefore, if $L \in J^{d^\prime} (Y)$ with
\begin{equation} \label{eq:ineq}
g_Y-g_X -1 -\frac{g_X^2}n \leq d^\prime \leq g_Y + g_X -1+\frac{g_X^2}n,
\end{equation}
then $\pi_* L$ is stable on $X$ and
$\pi_* (L \otimes \pi^*
M)$ is stable for any line bundle $M$ on $X$.  Thus,
for a generic line bundle $L_0\in J^d (Y)$ with
$d$ such that:

\noindent (i) $d>\operatorname{max}\{2g_Y-2+2n - (\operatorname{deg}
R_\pi)/2\ ,\ 2g_Y-2 \}$,

\noindent (ii) $d$ is equivalent (modulo $n$) to a number
$d^\prime$ such that
$d^\prime $
fulfills the condition  (\ref{eq:ineq}),

\noindent we have that $\pi _*(L_0)$ is stable and $\operatorname{deg} \pi_*(L_
0)=
\operatorname{deg} L_0-(\operatorname{deg} R_{\pi})/2>2ng_X$.
In the case $\operatorname{deg}\pi < \operatorname{max}\{ g_X(1+\sqrt
3)-1,2g_X+2\}$, it is
sufficient that
$\operatorname{deg} L_0$ fulfills
the condition
(i).
\end{em}
\end{rem}

\begin{thm}\label{th:beau}  If  $\pi_*(L_0)=E_0$ is stable and
$\operatorname{deg} (E_0) > 2ng_X$, then
$(\pi^*)^* (\cw_J)$ is $\Theta_X$-stable on $J(X)$.
\end{thm}

\begin{pf} For convenience of writing, we set up the following commutative
diagrams:
\begin{equation}\label{eq:diagram}
\xymatrix{ X & Y \ar[l]_{\pi} &  \\
 X\times J(X) \ar[u]_{p_X}
\ar[dr]_{p_{J(X)}}
 & Y\times J(X) \ar[u]^{q_1} \ar[l]_{\pi\times id} \ar[d]^{q_2}
\ar[r]^{id\times \pi^*} & Y\times J(Y)  \ar[d]_{p_2}
\ar[ul]_{p_1}  \\
    & J(X) \ar[r]^{\pi^*}  & J(Y)  }
\end{equation}
where the vertical morphisms are the natural projections.

Fix the Poincar\'e bundle  $\cp_X$ on $X\times J(X)$ normalized with
respect to the point $\pi(y_0)\in X$. Then,  the vector bundle
$(id\times\pi^*)^*\cp_Y$ on $Y\times J(X)$ is isomorphic to the bundle
$(\pi \times id)^* \cp_X$.
{}From diagram (\ref{eq:diagram}) and the base change
formula, we have
$$
\begin{aligned}
(\pi^*)^* (\cw_J) &\cong
q_{2*} \left( q_1^* L_0 \otimes (id\times \pi^*)^*
\cp_Y \right) \cong
p_{J(X)*} \left( (\pi\times id)_* \left(q^*_1 L_0\otimes (\pi\times
id)^*\cp_X\right) \right)\\
&\cong
p_{J(X)*} \left((\pi\times id)_* (q_1^* L_0) \otimes \cp_X \right)\cong  p_{J(X)*
}(p_X^*
(E_0)\otimes \cp_X).
\end{aligned}
$$

Since $E_0$ is stable of degree $d>2ng_X$, $p_X^* (E_0)\otimes \cp_X$ is a
family of stable
bundles parametrized by $J(X)$. Such family corresponds to an embedding of
the Jacobian in
the moduli space $\cm (n,d)$ of stable vector bundles of rank $n$ and
degree $d$. As in the
proof of Theorem 2.5 in
\cite{Li}, it follows that
$p_{J(X)*}(p_X^* (E_0)\otimes\cp_X)$ is $\Theta_X$-stable.
\end{pf}

With the hypothesis of the previous theorem, we obtain

\begin{cor}\label{cor:sta} The restriction $\cw _{\pi}$ of $\cw _J$ to $\pi
^*(J(X))$ is
$\Theta _{\pi}$-stable.
\end{cor}

\begin{pf}
The map $h:J(X)\rightarrow \pi^*(J(X))$ is an isogeny
such that $h^{-1}(\Theta _{\pi })\equiv n\Theta_X.$
Since $h^*(\cw _{\pi})\cong (\pi^*)^*(\cw _J)$,
we have from
Lemma 2.1 in \cite{bbn} that
$\cw _{\pi}$ is $\Theta _{\pi}$-stable.
\end{pf}

To study the restriction $\cw _P$ of $\cw _J$ to the Prym variety $P$
 we consider two cases, namely when the Abel-Prym map
$\widetilde{\rho}$ is not birational and when it is.
In both cases we can reduce our study to consider the
corresponding vector bundle over $(\widetilde{P},\varXi)$
where $f:\widetilde{P} \rightarrow P $ is a fixed isogeny as in Theorem
1.1. That is, if $j_P:P\hookrightarrow J(Y)$ is the natural
inclusion, then we will
denote by $\beta :\widetilde{P} \rightarrow J(Y) $
the composition map $f\circ j_P$ and by $\cw _{\widetilde{P}}$
the vector bundle $\beta^*(\cw _J)$. Actually,
$f^*(\cw _P) \cong \cw _{\widetilde{P}}.$

\begin{propn}\label{prop:simt} If $\cw_{\widetilde{P}}$
 is $\varXi$-stable (resp. simple),
then $\cw _P$ is $\Theta _P$-stable (resp. simple).
\end{propn}

\begin{pf} The stability again follows from
 Lemma 2.1 in \cite{bbn}. Suppose $\cw _{\widetilde{P}}$ is simple.
 Since  $\co_P \hookrightarrow f_* \co_{\widetilde P}$ is
injective, the map
$$
\begin{aligned}
\Hom(\cw_P , \cw_P) \hookrightarrow \Hom (\cw_P , \cw_P \otimes f_*
\co_{\widetilde P}) & \cong \Hom( \cw_P ,
f_* f^* \cw_P) \\ 
&\cong \Hom( \cw_{\widetilde P} ,\cw_{\widetilde P}) \cong {\Bbb C}
\end{aligned}
$$
is injective. 
Hence, $\Hom(\cw_P, \cw_P)\cong {\Bbb C}$ and $\cw_P$ is simple.
\end{pf}

Suppose $\widetilde{\rho}$ is not birational. As in section \S 1, let $Z$ be
the normalization of ${\widetilde{\rho}}(Y)$ and $\widetilde \pi: Y\to Z$ the
morphism induced by $\widetilde \rho$ of degree $\tilde n$.

\begin{propn}\label{prop:partsimp} If the Abel-Prym map $\widetilde \rho \colon
Y \to  \widetilde P$ is not birational onto
its image and
$\operatorname{deg} L_0> 2g_Y+4$, then the  Picard bundle $\cw_P$ is
$\Theta_P$-stable.
\end{propn}

\begin{pf} We consider separately each one of the cases where
 $\widetilde \rho$ is not birational (see Remark \ref{rem:iper}).

1) Suppose $n=2$.

By lemma \ref{lem:div},
$\operatorname{deg} \tilde \pi = 2$ and
$\widetilde{\rho}(Y)$ is of class one in $(\widetilde P, \varXi)$. Then,
by Matsusaka's Criterion, $Z=\widetilde{\rho}(Y)$ and
$(J(Z),\Theta_Z)\cong (\widetilde P, \varXi)$.

{}From the construction of the map $\widetilde \rho$, we have
$\operatorname{Nm} \tilde \pi = \tilde \mu$. Moreover, ${\tilde \mu}^t =
j_P \circ f = \beta$
(see section \S 1 in \cite{We}). Therefore the map $\tilde \pi ^* : J(Z)
\to J(Y) $
coincides with the map $\beta: \widetilde P \to J(Y)$ via the isomorphism
$J(Z)
\cong \widetilde P$. By definition of $\cw_{\widetilde{P}}$, we have
$$\cw _{\widetilde{P}}\cong (\widetilde{\pi}^*)^*(\cw _J).$$

In this case, $\operatorname{deg}\tilde \pi =2 <
\operatorname{max}\{ g_Z(1+\sqrt 3)-1,2g_Z+2\}$ because  $g_Z\geq 1$.
Therefore,
by Remark
\ref{rem:beau} and Theorem \ref{th:beau}, if
$\operatorname{deg} L_0 =d>2g_Y+2\geq\operatorname{max}
\{2g_Y+2-\operatorname{deg} (R_{\tilde \pi})/2,
2g_Y-2\}$, then
 $(\widetilde{\pi}^*)^*(\cw _J)\cong {\cw}_{\widetilde{P}}$
is $\varXi$-stable and from Proposition
\ref{prop:simt}, $\cw_P$ is
$\Theta_P$-stable.

2) Suppose $n=3$, $\tilde n = 3$, $g_Y = 4$, $g_X = 2$, $g_Z = 2$,
$R_\pi=R_{\tilde \pi}=0$.

Since $\widetilde{\rho}(Y)$ is of class one in $(\widetilde P, \varXi)$,
$Z=\widetilde{\rho}(Y)$ and
$(J(Z), \Theta_Z)\cong (\widetilde P, \varXi) $. In this case,
$\operatorname{deg}\tilde \pi =3 <
\operatorname{max}\{ g_Z(1+\sqrt 3)-1,2g_Z+2\}$. Hence, if
$\operatorname{deg} L_0
=d > 2g_Y +4=12$,
then $ {\cw}_{\widetilde{P}}$ is $\varXi$-stable and $\cw_P$ is
$\Theta_P$-stable. \par

3) Suppose $n=4$, $\tilde n= 2$, $g_Y=5$, $g_X=2$, $g_Z=3$,
$R_\pi = R_{\tilde
\pi}=0$.

We have the following commutative diagrams:
\begin{equation*}
 \xymatrix@=30pt{
 Y \ar @/^1.5pc/[rr]|{\,\widetilde \rho\,} \ar [r]^-{\widetilde \pi}\ar
@{_{(}->}[d]_-{\alpha_{y_{\text {o}}}}  &
 Z \ar[r] \ar @{_{(}->}[d]_-{\alpha_{\tilde {\pi}(y_{\text {o}})}} &
\widetilde {P}
\\
 J(Y) \ar[r]^-{\operatorname {Nm}\widetilde\pi} \ar @/_3.5pc/[rru] |{\,
\tilde \mu\, } &
 J(Z) \ar[ru]^-{u}   &
}
\end{equation*}
 where $u$ is the map induced by $Z\to \widetilde P$. Therefore,
$$\beta =j_P\circ f = \tilde \mu^t =
\tilde \pi^* \circ u^t\colon \widetilde P\longrightarrow
J(Y)$$
and $\operatorname{Im} {\tilde \mu}^t = P\subset \operatorname{Im} \tilde
\pi^*$.
Since $\dim P = 3 = \dim J(Z)= \dim \tilde \pi^* (J(Z))$, it follows that
$\tilde \pi^* (J(Z)) = P \subset J(Y)$.
Since $\operatorname{deg}\tilde \pi =2 <
\operatorname{max}\{ g_Z(1+\sqrt 3)-1,2g_Z+2\}$, if
$\operatorname{deg}L_0 =d >
2g_Y +2=12$ then, from Corollary \ref{cor:sta}, the restriction of $\cw _J$
to $\widetilde{\pi}^*(J(Z))$ is $\Theta_Y |_{\widetilde{\pi}^*(J(Z))}$-stable,
i.e. $W_P$ is $\Theta _P$-stable.
\end{pf}

\begin{rem} \begin{em} Observe that, in the previous proposition, the result 
holds for any line bundle $L_0$ (of degree $d> 2g_Y+4$) whereas in Remark 
\ref{rem:beau} the condition ``$\pi_*(L_0)$ is stable'' holds for a generic line 
bundle
$L_0$. If the Picard bundle $\cw_{L_0,P}$ corresponding to $L_0$ is stable, then 
$\cw_{L,
P}=p_{2*}(p_1^* L \otimes \cp_Y |_{Y\times P})$ is stable 
for any $L$ (of degree $d$). In fact, since $L\otimes L_0^{-1} \in J(Y)$, one  
can write
$L \cong L_0 \otimes \pi^* N \otimes M$ where $M\in P$  and $N\in J(X)$, 
therefore the two
Picard bundles are related as follows $\cw_{L,P}\cong \tau_M^* (\cw_{L_0\otimes \
pi^*
N,P})$, where $\tau_M\colon P\to P$ is the traslation by $M$. Since  $\pi_*(L_0)$ 
is
stable, $\pi_*(L_0\otimes \pi^*N)$ is 
stable. Hence $\cw_{L,P}$ is stable.
\end{em}\end{rem}

Before we consider the case when the map $\widetilde{\rho}$ is birational
we describe the bundle $\cw _{\widetilde{P}}$ in different way.

The abelian variety $(\widetilde P, \varXi)$ is principally polarized,
so we can identify
$\widetilde P$ with its dual abelian variety. The Poincar\'e bundle on 
$\widetilde P\times
\widetilde P$ is given by
$$\cq \cong m^* \co_{\widetilde P} (\varXi)\otimes q_1^*
\co_{\widetilde P} (\varXi)^\vee \otimes q_2^* \co_{\widetilde P}
(\varXi)^\vee ,$$ where $m$ is the multiplication law on
$\widetilde P $ and
$q_i$ the canonical projections of $\widetilde P \times
\widetilde P$ in the $i=1,2$ factor.

\begin{propn}\label{prop:funtrela2} The vector bundle
 $ {\cw}_{\widetilde{P}}$
is isomorphic to the bundle $q_{2*}
\left(q_1^*
\left({\widetilde
\rho}_* (L_0) \right)
\otimes \cq^\vee \right)$, where $q_i$ are the projection of
$\widetilde{P}\times \widetilde{P} $
in the $i=1,2$ factor.
\end{propn}

\begin{pf}Consider the commutative diagrams:
\begin{equation*}
\xymatrix{  & Y \ar[r]^{\widetilde{\rho}} & \widetilde{P} \\
Y\times P & Y\times \widetilde{P} \ar[l]_{id\times f}
\ar[r]^{\widetilde{\rho}\times id}
\ar[dr]_{p_2} \ar[u]^{p_1} &
\widetilde{P}\times \widetilde{P}\ar[d]^{q_2} \ar[u]_{q_1}\\
 &  & \widetilde{P} }
\end{equation*}

Let $\cp _{\widetilde{P}}$ be the line bundle $\cp _{\widetilde{P}}=(id\times 
(f\circ j_P))^*(\cp_Y)$ on $Y\times \widetilde P$, where $j_p:P\rightarrow J(X)$ 
is the natural
inclusion. From the normalization of $\cp_Y$, the restrictions of the line
bundles $\cp _{\widetilde{P}}$ and $({\widetilde \rho}\times id)^* \cq^\vee $  to 
$\{y_0\}\times\widetilde P$ are trivial and for
every $\tilde \xi \in \widetilde P$, the restrictions of the both bundles
to $Y\times \{\tilde
\xi\}$ are isomorphic to $f(\tilde \xi)$ by Lemma \ref{lem:lemcomp}. Therefore, 
$\cp
_{\widetilde{P}}\cong ({\widetilde \rho}\times id)^* \cq^\vee $. 

 From the projection formula and the base-change formula, we have
$$
\begin{aligned}
{\cw}_{\widetilde{P}} &\cong p_{2*} \left(p^*_1 L_0\otimes
{\cp _{\widetilde{P}}} \right)
\cong q_{2 *} \left( ({\widetilde \rho} \times id)_*
\left(p^*_1 L_0\otimes ( {\widetilde \rho}\times
id\right)^* \cq^\vee )\right) \\
 & \cong q_{2*}
\left( ( {\widetilde  \rho}\times id)_* (p_1^* L_0) \otimes \cq^\vee \right)
\cong q_{2*} \left(  q^*_1
\left( {\widetilde \rho}_{*} (L_0) \right)\otimes
\cq^\vee \right).
\end{aligned}$$
\end{pf}

\begin{propn} \label{prop:picsimp} If the Abel-Prym map $\widetilde {\rho}
$ is birational onto
its image then the restriction $\cw_P$ of the Picard bundle $\cw _J$ to the
Prym variety is simple.
\end{propn}

\begin{pf} From Proposition \ref{prop:simt} it is enough
to prove that ${\cw}_{\widetilde{P}}$ is simple.
The bundle $ {\cw}_{\widetilde{P}}$
is the Fourier-Mukai
transform of the sheaf ${\widetilde \rho}_* (L_0)$ by Proposition
\ref{prop:funtrela2}. Since
$\operatorname{deg}(L_0)>2g_Y-2$, the sheaf $q_1^* (\widetilde {\rho}_*
(L_0))\otimes \cq^\vee $ has only one non-null direct image. Then the
Fourier-Mukai transform is a complex
concentrated in degree zero and, therefore,
${\cw}_{\widetilde{P}}$ is simple if $\widetilde
\rho_* (L_0)$ is  simple (see \cite{Mu},
Corollary 2.5).

To prove that  $\widetilde \rho_* (L_0)$ is simple write
 $\widetilde \rho$ as a composite $j\circ \tilde \nu$ where
$\tilde \nu \colon Y\to \widetilde {\rho} (Y)$ is a birational morphism
between curves and
$j\colon {\widetilde \rho} (Y)\to \widetilde P$ is a closed embedding. Since
$\widetilde {\rho}$ is birational onto its image,  $\tilde \nu_* (L_0) $
is of rank $1$.
Therefore the sheaf $K$ defined by
$$0 \to K \to \tilde \nu^* \tilde \nu_* L_0 \to L_0 \to 0 $$
is a torsion sheaf. From the exact sequence
$$0\to \Hom (L_0,L_0)\to \Hom ( \tilde \nu^*\tilde \nu_* L_0, L_0)\to \Hom
(K, L_0).$$
and the fact that $ \Hom (K, L_0) = 0$, we obtain
$$\Hom(\tilde \nu_* L_0,\tilde \nu_* L_0)\cong \Hom ( \tilde \nu^*
\tilde \nu_* L_0, L_0) \cong
\Hom (L_0,L_0)\cong {\Bbb C}.$$

Since $j$ is a closed embedding, $j^* j_* E\cong E$ for any sheaf $E$
and applying
the adjunction formula again, we obtain that $j_* \tilde \nu_* (L_0)   \cong
\widetilde {\rho}_* L_0$
is simple.
\end{pf}

 From Proposition \ref{prop:partsimp}
 and Proposition \ref{prop:picsimp} we have that:

\begin{thm}\label{th:sim} The restriction $\cw_P$ of the
Picard bundle $\cw _J$ to the Prym variety  is
simple if $g_X\geq 2$ and
$\operatorname{deg}L_0>2g_Y+4$.
\end{thm}

\section{Picard bundles over moduli spaces}\label{sec:tres}

 We shall recall the construction and properties of spectral coverings
given by Beauville, Narasimhan and Ramanan in \cite{bnr}.

Let $X$ be a non-singular projective irreducible curve of genus
$g_X\geq 2$. Denote by  $\cm (n,d)$
the moduli space of stable vector bundles over $X$ of rank $n$ and degree
$d$. If $n$ and $d$ are coprime there exists a
universal family $\cu $ parametrized by $\cm (n,d)$ called the Poincar\'e
bundle. If $d> 2n(g_X-1)$ the Picard bundle $\cw$ is the direct image of
$\cu$ and it
is locally free. Denote by $\cw _{\xi}$ the restriction of $\cw$ to the
subvariety $\cm _{\xi}\subset \cm (n,d)$ determined by stable bundles with
fixed determinant
$\xi \in J^d(X)$.

Let $K$ be the canonical bundle over $X$ and $W= \bigoplus ^n_{i=1}H^0(X,K^i)$.
For every element $s=(s_1,\dots ,s_n)\in W$,
denote by $Y_s$ the associated spectral curve (see \cite{bnr}).
For a general $s\in W$, $Y_s$ is non-singular of genus $g_{Y_s} = n^2
(g_X-1)+1$ and
the morphism
$$\pi _s \colon Y_s \longrightarrow X$$
is of degree $n$.

In \cite{bnr} it is proved that if $\delta = d+n(n-1)(g_X-1)$ then there is an
open subvariety $T_{\delta}$ of  $J^{\delta}(Y_s)$ such
that the morphism $T_{\delta} \rightarrow \cm (n,d)$
defined by $L\mapsto \pi _{s*}(L)$
is dominant. Moreover, the direct image induces a dominant rational map
\begin{equation*}
\xymatrix{ f\colon P' \ar@{-->}[r] & \cm _{\xi} }
\end{equation*}
defined on an open subvariety $T'\subseteq P'$, where $P'$ is a translate
of the Prym
variety $P_s$ of $\pi _s$ (see \cite{bnr}, Proposition 5.7). The complement
of the open
subvariety $T'\subseteq P'$ is of codimension at least $2$. Actually,
$f\colon  T'
\longrightarrow {\cal M}_{\xi} $ is generically finite.

  Consider the following commutative diagram:
\begin{equation*}
\xymatrix{ Y_s\times T'\ar[r]^{\pi _s \times id} \ar[dr]_{q'_2} &
X \times T' \ar[r]^{id\times f} \ar[d]_{q_2} & X \times \cm
_{\xi}\ar[d]_{p_2}\\
 & T'\ar[r]^{f} & \cm _{\xi} }
\end{equation*}
where $q'_2,q_2,p_2$ are the projections to the second factor.

If ${\cal P}_{T'}$ is the restriction of the Poincar\'e bundle over
 $Y_s\times P^\prime $ to $Y_s\times T'$ and $\cu _{\xi}$ is the
restriction of the universal
bundle $\cu$ to $X \times \cm _{\xi}$, then, by the definition of $f$,
$$ (id \times f)^*(\cu _{\xi})\cong (\pi _s \times id)_*({\cal P}_{T'})\otimes
q_2^*(M)$$
for some line bundle $M$ over $T'$, which depends on the choice of $\cu$.
Therefore,
$$
f^*(\cw _{\xi})\cong (q_2)_*(id \times f)^*(\cu _{\xi})
\cong  (q_2)_*((\pi _s \times id)_*({\cal P}_{T'}))\otimes M
\cong (q'_2)_* ({\cal P}_{T'})\otimes M.
$$

Actually, $(q'_2)_* ({\cal P}_{T'})$ is just the restriction of
the Picard bundle $\cw _{P^\prime} $ over $P^\prime $ to $T'$.

\begin{thm}\label{th:mod} The Picard bundle $\cw _{\xi}$ on $\cm_\xi$ is
simple for
$d>n(n+1)(g_X-1)+6$ and $g_X\geq 2$.
\end{thm}

\begin{pf} Since $\operatorname{codim} (T^\prime) \geq 2$, we have
$$  \operatorname{End} (\cw_{P^\prime}) \cong H^0 (P', {\cal E}nd
(\cw_{P^\prime}))\cong H^0 (T', {\cal E}nd  (\cw_{P^\prime}|_{T'})) \cong
\operatorname{End} (f^* \cw_\xi) .$$

The Abel-Prym map for the spectral cover is birational onto its image if
$n\neq 2$ or $g_X>2$.  In this case, from Proposition \ref{prop:picsimp},
the bundle $\cw_{P'}$ is simple if
$\delta = d+n(n-1)(g_X-1)>2g_{Y_s}-2$.  Since the map $f\colon T'\to \cm_\xi$ is
 dominant and generically finite, as in the
 proof of Proposition \ref{prop:simt}, we deduce that
the Picard bundle $\cw _{\xi}$ on $\cm_\xi$ is
 simple for degree $d>n(n+1)(g_X-1)$.

 For $n=2$ and $g_X=2$, if the map $\widetilde \rho$ is not birational onto
its image, then
by the proof of Proposition \ref{prop:partsimp},
$\cw _{\xi}$ is simple for $d>n(n+1)(g_X-1) +6$. 
\end{pf}

\begin{rem}\label{rem:stabdes}\begin{em}
 Denote by $\Theta _{\xi}$ a generalized theta
divisor on $\cm _{\xi}$. By \cite{Li}, Theorem 4.3,
we have that
$$f^*({\cal O}(\Theta _{\xi}))\cong {\cal O}({\Theta _{P'}})|_{T'}$$
where $\Theta _{P'}$ is the restriction of the theta divisor of $J^\delta
(Y_s)$
to $P'$. Since $\cw_{P'}|_{T'} \cong f^* \cw_{\xi}$, by Lemma 2.1 in
\cite{bbn}, it
follows that if $\cw _{P'}$ is $\Theta _{P'}$-stable, then $\cw _{\xi}$ is
$\Theta _{\xi}$-stable. Moreover, from Theorem 2.2 in \cite{BHM},
$\cw_{P'}$ is
stable if $\widetilde \rho^*_L (\cw_{P'})$ is stable on $Y$ for a general
line bundle $L$.

\end{em}\end{rem}

Now we will focus on the case $n=2$ and $g_X=2$. $\cm(2,\xi)$ will denote 
the moduli space of stable rank $2$ vector bundles with fixed determinant 
$\xi$ (with $\operatorname{deg}\xi$ odd). We shall prove that it is 
possible to construct a spectral covering $\pi \colon Y_s \to X$ of degree $2$
 in such way that the curve $Y_s$ is hyperelliptic. Let $p\colon X \to 
{\Bbb P}^1$
 be the covering of degree $2$ given by 
the canonical bundle $K_X$. Let $s$ be a section of $H^0 (X, K_X^2)$ such that 
the spectral
covering $\pi\colon Y_s \to X$ of degree $2$ corresponding to $(0,s)$ 
is smooth, integral and such that the induced map from the Prym variety
associated to $Y_s\to X$ to the moduli space $\cm(2,\xi)$ is dominant. 
Observe that since $g_X=2$ we can write the section $s$ as a product $s=
s_1\cdot s_2$ where $s_i\in H^0(X, K_X)$ for $i=1,2$. If $P_i + Q_i$ is the
divisor associated to the section $s_i$, then $p(P_i) = p(Q_i) = a_i \in {\Bbb 
P}^1$. 

 Following the construction of cyclic coverings of curves given in
\cite{Go1},  the morphism $ \bar{X}=  {\Bbb P}^1\to {\Bbb P}^1$ of degree
$2$,  ramified  at the points ${a_1,a_2}$ is the cyclic covering defined by 
the construction data $ (D=a_1+a_2, L= {\mathcal O}_{{\Bbb P}^1}
(2))$. Let $C$ be the desingularization of the curve 
$X\times_{{\Bbb P}^1} \bar X$. From Theorem 2.13 in \cite{Go1} the map 
$C\to X$ is the cyclic covering associated to the data 
$(p^{-1} (D), p^*L) $, but since $p^{-1} (D) = \{P_1,Q_1,P_2,Q_2\}$ and
$p^*L \cong K_X^2$, then $Y_s \cong C$ and the map 
$Y_s\cong C \to \bar{X}={\Bbb P}^1$  is of degree $2$. Therefore the
curve $Y_s$ is hyperelliptic.

Let $P$ be the Prym variety associated to the degree $2$ covering 
$Y_s \to X$. Since the spectral curve $Y_s$ is hyperelliptic then the 
Abel-Prym morphism $\rho \colon Y_s \to P$ is not 
birational onto its image. As in the case 3) of the proof of Proposition 
\ref{prop:partsimp}
its follows that, in this case, if $\delta = d+n(n-1)(g_X-1)> 2g_{Y_s}+2=12$,
then the Picard bundle on the Prym variety $P$ is stable with respect to the 
restriction of the theta divisor. From Remark \ref{rem:stabdes} the stability 
of  Picard bundles on $\cm_X (2,\xi)$ is ensured. 
 
\begin{thm} Let $X$ be a smooth projective irreducible 
curve of genus $2$ and let $\cm_X (2,\xi)$ be the moduli space of 
of rank $2$ stable vector bundles on $X$ with fixed determinant $\xi$. 
If $d= \operatorname{deg} (\xi) > 10$ and $d$ odd, then the Picard
bundle $\cw_{\xi}$ on $\cm_X (2,\xi)$ is stable with respect to  
the theta divisor.
\end{thm}

\begin{rem}\begin{em}
After this paper was finished, I. Biswas and T. G\'omez informed us that they 
have obtained the stability of the Picard bundle in the case of the 
moduli space $\cm_X (2, \xi)$ of rank $2$ stable bundles over a smooth curve 
$X$ of genus $g_X\geq 3$ such that $d=\operatorname{deg} \xi \geq 4g_X-3$ 
and odd (\cite{BG}).
\end{em}\end{rem}


\end{document}